\documentclass[12pt,leqno]{article}
\usepackage{amsmath, amsthm, amsfonts, amssymb, color}
\usepackage{mathrsfs}
\newtheorem{thm}{Theorem}[section]

\newtheorem{lem}[thm]{Lemma}
\newtheorem{prp}[thm]{Proposition}

\theoremstyle{definition}

\newcommand{\scr}[1]{\mathscr #1}
\definecolor{wco}{rgb}{0.5,0.2,0.3}

\numberwithin{equation}{section} \theoremstyle{remark}

\newcommand{\ua}{\uparrow}

\title{{\bf Transportation Cost Inequality on  Path Spaces with Uniform Distance}
\footnote{Supported in part by
NNSFC(10121101) and RFDF(20040027009)}}

\author{{\bf Shizan Fang$^{a,b}$, Feng-Yu Wang$^{a,c}$\footnote{wangfy@bnu.edu.cn}},\ \ {\bf Bo Wu$^{a}$}\\
\footnotesize{$^a$ School  of Mathematical Science \&  Lab. Math. Com. Sys.,}\\
 \footnotesize{ Beijing Normal
University, Beijing 100875, China}\\
{\footnotesize $^b$ I.M.B. B.P. 47870, Universit\'e de Bourgogne, Dijon, France}\\
\footnotesize {$^c$ WIMCS, University of
           Wales Swansea, Singleton Park, Swansea, SA2 8 PP, UK }}
\begin{document}
\maketitle

\def\R{\mathbb R}  \def\ff{\frac} \def\ss{\sqrt} \def\BB{\mathbb
B}
\def\N{\mathbb N} \def\kk{\kappa} \def\m{{\bf m}} \def\C{\scr C}
\def\dd{\delta} \def\DD{\Delta} \def\vv{\varepsilon} \def\rr{\rho}
\def\<{\langle} \def\>{\rangle} \def\GG{\Gamma} \def\gg{\gamma}
  \def\nn{\nabla} \def\pp{\partial} \def\tt{\tilde}
\def\d{\text{\rm{d}}} \def\bb{\beta} \def\aa{\alpha} \def\D{\scr D}
\def\E{\mathbb E} \def\si{\sigma} \def\ess{\text{\rm{ess}}}
\def\beg{\begin} \def\beq{\begin{equation}}  \def\F{\scr F}
\def\Ric{\text{\rm{Ric}}} \def\Hess{\text{\rm{Hess}}}\def\B{\mathbb B}
\def\e{\text{\rm{e}}} \def\ua{\underline a} \def\OO{\Omega} \def\sE{\scr E}
\def\oo{\omega}     \def\tt{\tilde} \def\Ric{\text{\rm{Ric}}}
\def\cut{\text{\rm{cut}}} \def\P{\mathbb P} \def\ifn{I_n(f^{\bigotimes n})}
\def\C{\scr C}      \def\aaa{\mathbf{r}}     \def\r{r}
\def\gap{\text{\rm{gap}}} \def\prr{\pi_{{\bf m},\varrho}}  \def\r{\mathbf r}
\def\Z{\mathbb Z} \def\vrr{\varrho} \def\ll{\lambda}
\def\ra{\mapsto} \def\MATHBF{\mathbf}
\def\L{\scr L}\def\Tt{\tt} \def\TT{\tt}
\def\Ric{{\rm Ric}} \def\LL{\Lambda}
\def\ric{{\rm ric}}
\def\FC{{\mathcal FC}_b^\infty}
\def\E{{\mathbb{E}}}
\def\R{\mathbb{R}}
\def\H{\mathbb{H}}
\def\Id{{\rm Id}}
\def\Hess{{\rm Hess}}

\begin{abstract} Let $M$ be a complete Riemnnian manifold and $\mu$ the distribution of the diffusion
process generated by $\ff 1 2\DD+Z$ where $Z$ is a $C^1$-vector
field. When $\Ric-\nn Z$ is bounded below and $Z$ has, for
instance, linear growth, the transportation-cost inequality with
respect to the uniform distance is established for $\mu$ on the
path space over $M$. A simple example is given to show the
optimality of the condition.
\end{abstract}
\



\vskip 2cm
\section{Introduction}

Since Talagrand \cite{T} found his transportation cost inequality
for the Gaussian measure on $\R^d$, this inequality have been
established on finite- and infinite-dimensional spaces with respect
to many different distances (i.e. cost-functions); see \cite{V} for
historical comments and  references. For instance, on the path space
of a diffusion process on a complete Riemannian manifold, the
inequality holds with respect to the $L^2$-distance provided the
curvature of the diffusion is bounded below, and it holds with
respect to the intrinsic distance induced by the Malliavian gradient
provided the curvature is bounded (see \cite{Wa1, Wa2}). See also
\cite{FU, Wu1, Wu2} for the study of  diffusion path spaces over
$\R^d$, and \cite{FSh, Sh} for the study on path and loop groups.
The purpose of this paper is to search for a reasonable curvature
condition   such that the Talagrand inequality holds for the
distribution of the corresponding diffusion process with respect to
the uniform distance on the path space.

Let $(M,g)$ be a connected complete Riemannian manifold of
dimension $d$. Consider the diffusion operator
$L=\frac{1}{2}\Delta_M +Z$ for a $C^1$-vector field $Z$.  Assume
that
\begin{equation}\label{4.1}
\Ric - \nabla Z\geq -K.
\end{equation}
Let $o\in M$ and $T>0$ be fixed. Let  $H: TM\to TO(M)$  the
horizontal lift, where $O(M)$ is the orthonormal frame bundle over
$M$. Consider the stochastic differential equation on $O(M)$:

\begin{equation}\label{4.2}
\d u_t(w)=\sum_{i=1}^d H_i(u_t(w))\circ \d w_t^i + \ff 1 2
H_Z(u_t(w))dt,\ u_0\in \pi^{-1}(o),
\end{equation} where $w_t =(w_t^i: 1\le i\le d)$ is the Brownian motion on
$\R^d$ and $H_i(u):= H_{u e_i}, 1\le i\le d.$ Here and in what
follows, $\{e_i\}_{i=1}^d$ is the canonical orthonormal basis on
$\R^d$. Let

\beg{equation*}\beg{split} &W_0(\R^d):=\{w\in C([0,T]; \R^d):\
w_0=0\},\\
&\mathbb H:= \bigg\{h\in W_0(\R^d):\ \|h\|_\H^2:= \int_0^T\dot
h_s^2\d s<\infty\bigg\},\end{split}\end{equation*} and $\P$ is the
standard Wiener measure. Let $\pi: O(M)\to M$ be the canonical
projection. Then $\gamma_t(w):=\pi u_t(w)$ is the $L$-diffusion
process on $M$ starting from $o$, which is non-explosive under the
condition (\ref{4.1}). Let $\mu$ be the law of $w\ra \gamma(w)\in
W_o(M):= \{\gg\in C([0,T]; M):\ \gg_0=o\}.$

Let $\rr$ be the Riemannian distance on $M$.  For two probability
measures $\mu_1,\mu_2$ on $W_o(M)$, let
$W_{2,d_\infty}^2(\mu_1,\mu_2)$ be the $L^2$-Wasserstein distance
between them induced by the uniform norm

$$d_\infty(\gg,\eta):=\sup_{t\in [0,T]}\rr(\gg_t,\eta_t),\ \ \ \gg,\eta\in
W_o(M).$$    More precisely,

\begin{equation}\label{(A.5)}
W_{2,d_\infty}^2(\mu_1,\mu_2):=\inf_{\hat\pi\in
\C(\mu_1,\mu_2)}\int_{W_o(M)\times
W_o(M)}d_\infty^2(\gamma,\eta)\hat\pi(d\gg,d\eta)
\end{equation}
where $\C(\mu_1,\mu_2)$ is the set of all probability measures on
$W_o(M)\times W_o(M)$ with marginal distributions $\mu_1$ and
$\mu_2$. The main result of the paper is the following:

\beg{thm}\label{T1.1}  Assume  $(\ref{4.1})$ hold for some $K\ge
0$ and let  $\rr_o=\rr(o,\cdot).$ If $ |Z|\le \psi\circ \rr_o$ for
some strictly positive and increasing $\psi\in
C^\infty([1,\infty))$ with

$$\int_0^\infty \ff 1 {\psi(s)}\d s =\infty,$$ then \beq\label{W4} W_{2,d_\infty}^2(F\mu,\mu)\leq
2\frac{\e^{KT}-1}{K}\mu(F\log{F}),\ \ \ F\ge 0, \mu(F)=1.
\end{equation} \end{thm}

\

To prove this result, we could start from the log-Sobolev inequality
for damped gradients $\tilde D$ (\ref{4.8}) below. To this end, one
would like to follow the line of \cite{BGL} by studying the
Hamilton-Jacobi semigroup $Q_t$ induced by the uniform norm
$d_\infty$:

$$(Q_tF)(\gamma)=\inf_{\eta\in W_o(M)}\Bigl\{F(\eta)+
\frac{1}{2t}d_\infty^2(\gamma,\eta)\Bigr\}.$$ By \cite{Sh}, $Q_t$
preserves the class of $d_\infty$-Lipschitz functions. Then,
according to the argument of \cite{BGL}, to derive the desired
transportation cost inequality from the log-Sobolev inequality
(\ref{4.8}), it remains to prove that

$$D_t^+ Q_t F:=\limsup_{s\downarrow 0} \ff{Q_{t+s}F-Q_tF}s\le - C \int_0^T |\tt D_s Q_tF|^2\d s$$
for some constant $C>0$, the inequality for which we are actually
in position to prove if $Z=0$.

So, in this paper we shall follow the line of \cite{Wa2} using
finite-dimensional approximations. To make the corresponding finite
approximate metric continuous, we have to first assume that the
Ricci curvature is $C_b^1$, the curvature tensor $\Omega$ is $C_b^0$
and the drift is $C_b^2$. So, to finish the prove,   we adopt one
more approximation argument on the Riemannian metric and the drift
to fit the above regularity assumption. To realize the second
approximation procedure we need the growth condition of $|Z|$ stated
in Theorem \ref{T1.1}. On the other hand, however, since the growth
of $|Z|$ is not included in the inequality (\ref{W4}), we believe
that it is  technical rather than necessary.

To conclude this section, we present below a simple example to show
that the condition in Theorem \ref{T1.1} for (\ref{W4}) is sharp.

\paragraph{Example 1.1.} Let $M=\R^d$ and $Z=\nn V$ for

$$V(x):= (1+|x|^2)^\dd,\ \ \ x\in \R^d,$$ where
$\dd\ge 0$ is a constant. Let $T>0$ and $o=0\in\R^d$ be fixed. We
claim that there exists a constant $C>0$ such that

\beq\label{W5} W_{2,d_\infty}^2(F\mu,\mu)\leq C\mu(F\log{F}),\ \ \
F\ge 0, \mu(F)=1
\end{equation} holds if and only if either $\dd\le 1.$
Indeed, for $\dd\le 1$  $\Ric-\nn Z=-\Hess_V$ is bounded from
below and $|Z|$  has   at most linear growth. So, (\ref{W5})
follows from Theorem \ref{T1.1}. On the other hand, it is
well-known that (\ref{W5}) implies

\beq\label{W6}\mathbb E \exp\Big[\ll\sup_{t\in
[0,T]}|\gg_t|^2\Big]<\infty\end{equation} holds for some $\ll>0,$
where $\gg_t$ is the $L$-diffusion process starting from $0$.
Indeed, according to \cite{Wu1}, this concentration property is
equivalent to the weaker $L^1$ transportation cost inequality:

$$W_{1,d_\infty}^2(F\mu,\mu)\leq C\mu(F\log{F}),\ \ \
F\ge 0, \mu(F)=1$$ for some constant $C>0$, where

$$W_{1,d_\infty}(\mu_1,\mu_2):=\inf_{\hat\pi\in
\C(\mu_1,\mu_2)}\int_{W_o(M)\times
W_o(M)}d_\infty(\gamma,\eta)\hat\pi(d\gg,d\eta)\le
W_{2,d_\infty}(\mu_1,\mu_2).$$ It is easy well-known that if
$\dd>1$ then the diffusion process is explosive so that (\ref{W6})
does not hold for any given $\ll>0.$

\

The remainder of the paper is organized as follows. In Section 2
we prove Theorem \ref{T1.1} under an additional assumption on
bounded geometry (see $(H)$ below), which in particular implies
the regularity of finite-dimensional metrics induced by
conditional expectations of the damped gradient. For readers'
convenience to follow the main points of the proof, we address the
proof of this regularity property in the Appendix at the end of
the paper. Then a complete proof of Theorem \ref{T1.1} is
presented in Section 3 by constructing  Riemannian manifolds
$\{(M_n,g_n): n\ge 1\}$ and operators $\{L_n: n\ge 1\}$, which
satisfy the assumption $(H)$ and approximate the original
Riemannian manifold and $L$ in a good way. Since the intrinsic
distance induced by the damped gradient on $W_o(M)$ is heavily
dependent of the geometry of $M$, it is not consistent through our
approximation. Finally, in Section 4 we extend Theorem \ref{T1.1}
to the free path space.

\section{The case with bounded geometry}

  In this section shall
assume that

\ \newline $(H)$\ \ \ $\Ric \in C_b^1,\ \Omega\in C_b^0 \
\hbox{and}\ Z\in C_b^2.$

\ \newline It is known that under $(H)$ the measure $\mu$ is
equivalent to the Wiener measure (see \cite{CHL}). It is also known
that the filtration generated by $\{u_s(w);\ s\leq t\}$ coincides
with the one generated by $\{\gamma_s(w); s\leq t\}$; they are both
equal to the natural filtration ${\mathcal N}_t$ generated by
$\{w_s;\ s\leq t\}$ (see \cite{Dr,TW}). For $F\in\FC,$ i.e.
\begin{equation}\label{1.4}
 F(\gamma)=f(\gamma_{s_1},\cdots, \gamma_{s_N}),\ 0<s_1<
\cdots< s_N\leq T,
\end{equation}
for some $N\ge 1$ and  $f\in C_b^\infty(M^N),$ we  define

$$(D_sF)(\gamma)=\sum_{j=1}^N u_{s_j}^{-1}(\partial_j f){\bf
1}_{\{s<s_j\}},$$ where $\pp_j$ is the gradient with respect to the
$j$-th component.

Throughout the paper, for any $p$-tensor $\scr T$ on $M$, let $\scr
T^\#: O(M)\to \scr L(\R^{d\times p}; \R)$ with

$$\scr T^\#(u)(a_1,\cdots, a_p)= \scr T(u a_1,\cdots, u a_p),\ \
\ \ a_1,\cdots, a_p\in \R^d, u\in O(M).$$

Now, let $\Ric_Z= \Ric -\nn Z$ and $\Ric_Z^\#$ be defined   for
$\scr T= \Ric_Z.$ Consider the following resolvent equation on $\scr
L(\R^d;\R^d)$:

\begin{equation}\label{4.4}
\frac{\d Q_{t,s}}{\d t}=-\frac{1}{2}\Ric_Z^\#(u_t) \,Q_{t,s},\quad
t\geq s>0;\quad Q_{s,s}=\Id.
\end{equation} By (\ref{4.1}),

\beq\label{A*} \|Q_{t-s}\|\le \e^{K(t-s)/2},\ \ \ t\ge
s>0,\end{equation} where $\|\cdot\|$ is the operator norm on $\R^d$.
Following \cite{FM}, we define the damped gradient

\begin{equation}\label{2.6}
(\tilde D_sF)(\gamma)=\sum_{j=1}^N Q_{s_j,s}^*
u_{s_j}^{-1}(\partial_j f){\bf 1}_{\{s<s_j\}},\ \ \
F(\gg)=f(\gg_{s_1},\cdots, \gg_{s_N}),
\end{equation} where $Q_{s_j,s}^*$ is the adjoint  of $Q_{s_j,s}.$
Then there holds the following integration by parts
formula

\beq\label{parts}  \E\Bigl(\int_0^T \bigl<\tilde D_sF,\dot
h(s)\bigr>\, \d s\Bigr) =\E\bigg(F\int_0^T \bigl<\dot h(s), \d
w_s\bigr>\bigg),\ \ \ h\in \mathbb H.
\end{equation} Indeed, letting
$\tilde h$ solve

\begin{equation}\label{2.4}
\dot{\tilde h}(t)+\frac{1}{2}\Ric_Z^\#(u_t)\,\tilde h(t)=\dot h(t),\
\ \tilde h(0)=0,
\end{equation} we have

\begin{equation*}
\begin{split}
\int_0^T \bigl<D_sF,\dot{\tilde h}(s)\bigr>\, \d s & =\sum_{j=1}^N
\bigl<u_{s_j}^{-1}(\partial_j f), \tilde h(s_j)\bigr>=\sum_{i=1}^N
\int_0^{s_j}
\Big\<u_{s_j}^{-1}(\pp_j f), \ff{\d}{\d s} Q_{s_j,s}\tt h(s)\Big\>\d s\\
&\hskip -10mm=\sum_{j=1}^N \int_0^{s_j}
\bigl<u_{s_j}^{-1}(\partial_j f), Q_{s_j,s}\dot h(s)\bigr>\, \d s\\
&\hskip -10mm=\sum_{j=1}^N \int_0^T
\bigl<Q_{s_j,s}^*u_{s_j}^{-1}(\partial_j f){\bf 1}_{\{s<s_j\}}, \
\dot h(s)\bigr>\, \d s =\int_0^T \<\tt D_s F, \dot h(s)\>\d s.
\end{split}
\end{equation*} Then
(\ref{parts}) follows from the known integration by parts formula
for the Malliavian gradient (see \cite{Bi, Dr, FM}). Under the
hypothesis $(H)$, we can use the approach \cite{CHL} to get
following log-Sobolev inequality (see also \cite{ELJL} for a
possible degenerate diffusion)

\begin{equation}\label{4.8}
 \mu(F^2\log F^2)\leq 2
\mu\Bigl(\int_0^T |\tilde D_sF|^2 \d s\Bigr),\ \ F\in\FC,
\mu(F^2)=1.
\end{equation}
Indeed, under our notations the last formula on page 75 of
\cite{CHL} (see Section 3 therein for the case with drift) becomes

$$H_t = \mathbb E\bigg( D_t F -\ff 1 2 \int_t^T Q_{s,t}^*
\Ric_Z^\#(u_s) D_s F \d s\bigg| \mathcal N_t\bigg)=\mathbb E(\tilde
D_t F|\mathcal N_t),$$ where the last equation follows from the
above relationship between the gradient and the damped gradient.
Then, replacing $F$ by $F^2$ in the second formula on page 75 in
\cite{CHL} and noting that $\mathbb E F^2=1$, we obtain

$$\mathbb E F^2\log F^2 \le 2\mathbb E\int_0^T \ff{[\mathbb
E(F\tilde D_t F |\mathcal N_t)]^2}{\mathbb E(F^2|\mathcal N_t)}\d
t\le 2 \mathbb E \int_0^T |\tilde D_t F|^2 \d t,$$ which is nothing
but (\ref{4.8}).

\

We shall derive the desired transportation-cost inequality from this
log-Sobolev inequality. It was observed by \cite{OV} (see also
\cite{BGL, Wa2}) that the log-Sobolev inequality  on a
finite-dimensional manifold implies the corresponding
transportation-cost inequality with respect to the intrinsic
distance of the associated Dirichlet form, which has been recently
extended in \cite{Sh} to an abstract setting under certain
assumption on the corresponding Hamilton-Jacobi semigroup. Since
this assumption does not directly apply to our present situation, we
shall adopt an approximation argument as in \cite{Wa2}. To this end,
we first reduce (\ref{4.8}) to a finite-dimensional log-Sobolev
inequality, which  implies a finite-dimensional transportation-cost
inequality; then  pass to the infinite-dimensional setting by taking
limit with respect to a sequence of partitions of $[0,T]$. Note that
the role of (\ref{4.8}) is only intermediate here, used throughout
the bounded geometry approximation; the constants involved will be
well behaved when the uniform distance will be taken into account.

\subsection{The finite-dimensional setting}

Let $I=\{0<s_1\cdots <s_N\leq T\}$ be a partition of $[0,T]$. Let

$$ \LL_I(\gg):= (\gamma(s_1), \cdots,
\gamma(s_N)),\ \ \ \gg\in W_o(M)$$ be the projection from $W_o(M)$
onto the product manifold $M^I.$ Then $\mu_I:=(\Lambda_I)_*\mu$ has
a   smooth and strictly positive density with respect to the
Riemannian volume $dx_1\cdots dx_N$ on $M^I$. For
$F=f\circ\LL_I\in\FC,$  it follows from (\ref{2.6}) that

\begin{equation}\label{A1}
\int_0^T|\tilde D_sF|^2\,\d s =\sum_{i,j=1}^N \int_0^{s_i\wedge
s_j}\<u_{s_j}Q_{s_j,s}Q_{s_i,s}^* u_{s_i}^{-1}\,\partial_if,\
\partial_j f\>_g\d s.
\end{equation}
Let $M^I\ni z\mapsto \P(z,\cdot)$ be the regular conditional
distributions of $\P$ given $\LL_I\circ\gg$. We define the linear
operator $A^I(z)$ on $T_z M^I$ by

\begin{equation}\label{A.7}
\<A^I(z)X, Y\>_{g^I} =\int_{W_0(\R^d)} \Bigl(\sum_{i,j=1}^N
\int_0^{s_i\wedge
s_j}\bigl<u_{s_j}Q_{s_j,s}Q_{s_i,s}^*u_{s_i}^{-1}X_i(z ) ,\ Y_j(z)
\bigr>_g\, \d s \Bigr) \d \P(z, \cdot)
\end{equation} for $X,Y\in T_z M^T$, where $g^I$ is the product Riemannian metric on $M^I$ and
 $X_i$ and $Y_j$ are the $i$-th and $j$-th
components of $X$ and $Y$ respectively. By Propositions \ref{lem4.3}
and \ref{lem4.4} below, $A^I$ is uniformly   positive definite and
has a continuous version, denoted again by $A^I$. Moreover, since  a
continuous mapping on $T M^I$ can be uniformly approximated by
smooth ones, in the sequel we may and do assume that $A^I$ is
smooth.

Noting that for   $F= f\circ \LL_I\in \FC$ we have

\begin{equation*}
\E\Bigl(\int_0^T|\tilde D_sF|^2\,\d s \Bigr)=\int_{M^I}\<A^I \nn^I
f, \nn^I f\>_{g^I} \d\mu_I,
\end{equation*} where $\nn^I$ is the gradient operator induced by $g^I$ on $M^I$, it follows from (\ref{4.8}) that

\begin{equation}\label{log}
\mu_I(f^2\log f^2) \le 2 \mu_I(\<A^I\nn^I f, \nn^I f\>_{g^I}),\ \ \
f\in C_b^\infty(M^I), \mu_I(f^2)=1.
\end{equation}
Now, let $\rho_I$ be the Riemannian distance induced by $A^I$ on
$M^I$. We have
\begin{equation}\label{4.16}
\rho_I(z,z')=\sup\bigl\{|f(z)-f(z')|:\  f\in C_b^1(M^I),
\bigl<A^I\nabla^If,\nabla^If\bigr>_{g^I}\leq 1\bigr\}.
\end{equation}
Since $g$ is complete, $(H)$ and Proposition \ref{lem4.4} below
imply the completeness of $\rr_I$. Therefore, by \cite[Theorem
1.1]{Wa1} with $p=2$ (see also \cite{OV, BGL}), (\ref{log}) implies
\begin{equation}\label{4.17}
W_{2,\rho_I}^2(f\mu_I,\mu_I)\leq 2\mu_I(f\log{f}),\quad f\geq 0,
\mu_I(f)=1.
\end{equation}
We are now ready to prove the main result of the paper under the
assumption $(H).$

\beg{prp}\label{P1} Assume $(\ref{4.1})$ and  $(H)$. Let $d_I(z,
z'):= \max_{1\le i\le N} \rr(z_{s_i}, z'_{s_i}),\ z,z'\in M^I.$ We
have

\beq\label{W3} W_{2,d_I}^2(f\mu_I,\mu_I)\leq
2\frac{\e^{KT}-1}{K}\mu_I(f\log{f}),\ \ \ f\ge 0, \mu_I(f)=1.
\end{equation} \end{prp}

\beg{proof} By (\ref{4.17}), we it suffices to prove that
\begin{equation}\label{4.19}
d_I^2\leq \frac{\e^{KT}-1}{K}\, \rho_I^2.
\end{equation} Obviously,

\beq\label{A2} d_I(z,z')= \sup\Big\{|f(z)-f(z')|:\ f\in
C_b^\infty(M^I), \sum_{j=1}^N |\pp_j f|_g\le 1\Big\},\ \ \ z,z'\in
M^I.\end{equation} Next, by (\ref{A1}) and the definition of $A^I$,
we have

\beq\label{A**}\rr_I(z,z')\ge \sup\bigg\{|f(z)-f(z')|:\ f\in
C_b^\infty(M^I), \int_0^T |\tt D_s F|^2\d s\le
1\bigg\}\end{equation} for $z,z'\in M^I.$ Finally, for $f\in
C_b^\infty(M^I)$ and $F= f\circ \LL_I$, (\ref{A1}) and (\ref{A*})
imply

$$\int_0^T |\tt D_s F|^2\d s \le \ff{\e^{KT}-1}{K} \sum_{i,j=1}^N
|\pp_i f|_g |\pp_j f|_g = \ff{\e^{KT}-1}K \Big(\sum_{j=1}^N |\pp_j
f|_g\Big)^2.$$  Therefore, (\ref{4.19}) follows from (\ref{A2}) and
(\ref{A**}).
\end{proof}

\subsection{The infinite-dimensional case}

\beg{prp}\label{P2} Assume $(H)$. Then $(\ref{4.1})$ implies
$(\ref{W4})$.  \end{prp}

\beg{proof} Since $\FC$ is dense in $L^1(\mu)$, it suffices to prove
(\ref{W4}) for nonnegative $F\in \FC$ with $\mu(F)=1.$ Let $F=
f\circ \LL_I$ for some partition $I$ of $[0,T]$ and nonnegative
$f\in C_b^\infty(M^I)$ with $\mu_I(f)=1.$ Take a sequence of
partitions $\{I_n\}$ of $[0,T]$, which is finer and finer such that
$I_n\supset I$ for all $n\ge 1$ and  $\cup_{n\geq 1}I_n$ is dense in
$[0,T]$. Let $$\tt d_{I_n}(\gg,\eta)= d_{I_n}(\LL_{I_n}(\gg),
\LL_{I_n}(\eta)),\ \ \ \gg,\eta\in W_o(M).$$ Then

\begin{equation}\label{4.21}
\tt d_{I_n}\uparrow d_\infty\ \ \text{as}\  \ n\uparrow\infty.
\end{equation}

Since $I_n\supset I$, we may regard $f$ as a function on $M^{I_n}$
depending only on components in $M^I$ so that $\mu_{I_n}(f)=1$ and
$\mu_{I_n}(f\log{f})=\mu(F\log{F})$ for all $n\geq 1$. By
(\ref{W3}), for any $n\geq 1$, there exists a coupling measure
$\tilde\pi_n\in C(f\mu_{I_n},\mu_{I_n})$ such that (cf. \cite{Ra})
\begin{equation}\label{4.22}
\begin{split}
\int_{M^{I_n}\times M^{I_n}}d_{I_n}^2 \d\tilde\pi_n \leq
2\ff{\e^{KT}-1}K \mu(F\log{F}).
\end{split}
\end{equation}
For any $n\ge 1$, let $\mu(z, \cdot)$ (respectively
$(F\mu)(z,\cdot)$) be the regular conditional distribution of $\mu$
(respectively $F\mu$) given $\LL_{I_n}=z \in M^{I_n}$.  Then,
according to \cite[page 187]{Wa2} (see also \cite[page 353]{FW}),

$$\hat \pi_n (\d\gg, \d\eta):=\int_{M^{I_n}\times M^{I_n}}
(F\mu)(z; \d\gg) \mu(z',\d\eta)\tt \pi_n(\d z, \d z') \in \scr
C(F\mu,\mu),\ \ \ n\ge 1.$$ Moreover, it is easy to see that
(\ref{4.22}) implies

\beq\label{W10} \int_{W_o(M)\times W_o(M)} \tt d_{I_n}^2 \d
\hat\pi_n \le 2\ff{\e^{KT}-1}K \mu(F\log{F}).
\end{equation} Since as  explained on page 187 of \cite{Wa2} that
$\scr C(F\mu,\mu)$ is tight and closed under the weak topology,  up
to a  subsequence   $\hat \pi_n\to \hat \pi$ weakly for some $\hat
\pi\in \scr C(F\mu,\mu)$ as $n\to\infty.$  Then, for any $N>0$, it
follows from (\ref{W10}) and the monotonicity of $\tt d_{I_n}$ in
$n$ that

\beg{equation*}\beg{split}&\int_{W_o(M)\times W_o(M)} (\tt
d_{I_N}^2\land N) \d\hat\pi = \lim_{n\to\infty} \int_{W_o(M)\times
W_o(M)} (\tt d_{I_N}^2\land N) \d \hat\pi_n\\
&\le \int_{W_o(M)\times W_o(M)} \tt d_{I_n}^2 \d \hat\pi_n\leq
2\ff{\e^{KT}-1}K \mu(F\log{F}).\end{split}\end{equation*} Therefore,
the proof is finished by taking $N\to\infty$ and using
(\ref{4.21}).\end{proof}

\section{Proof of Theorem \ref{T1.1}}

To prove Theorem \ref{T1.1} from Proposition \ref{P2}, we shall
constructed a sequence of metrics $\{g_n\}$ and drifts $\{Z_n\}$
satisfying $(H)$ and $\Ric_n -\nn_n Z_n\ge -K_n$ with $K_n\to K$
and $\mu_n\to\mu$, where $\Ric_n, \nn_n$ are the Ricci curvature
and the Levi-Civita connection induced by $g_n$, and $\mu_n$ is
the distribution of the diffusion process generated by $L_n:= \ff
1 2 \DD_n+Z_n$. Here, we will take  $g_n$ as  conformal changes of
$g$. So, we first study the conformal change of metric.

\subsection{Conformal changes of metric for $(H)$}

In this subset, we prove that the conformal change of metric used in
\cite{TW} satisfies the assumption $(H)$. More precisely, let $f\in
C_0^\infty(M)$ with $0\le f\le 1$ such that $M':= \{f>0\}$ is a
non-empty open set. Then, according to \cite{TW}, $M'$ is a complete
Riemannian manifold under the metric $g':= f^{-2}g$, and

$$L':= f^2 L=\ff 1 2 \DD' +Z'$$
for $Z'= f^2 Z+ \ff {d-2} 2 \nn f^2,$ where $\DD'$ is the Lapalcian
  induced by $g'.$ Let $\scr X(M')$ be the set of all smooth vector fields on $M'$,
  and $\scr X_b^p(M',g')$ (respectively, $\scr X_b^p(M',g)$) the set of all $C_b^p$ vector
  fields on $M'$ with respect to the metric $g'$ (respectively,
  $g$). Moreover,
Let$\nn'$ be the Levi-Civita connection
  on $(M',g').$
We have (see \cite[Theorem 1.159 a)]{E}) \beq\label{N*1} \nn_X' Y=
\nn_X Y- \<X,\log f\>_g X-\<Y, \log f\>_g Y +\<X,Y\>_g \nn \log
f,\ \ \ X, Y\in \scr X(M').\end{equation} Moreover, letting
$\Ric'$ be the Ricci curvature
 on $(M',g')$, by \cite[Theorem 1.159 d)]{E} we have (note that the
Laplacian therein equals to our  $-\DD$)

\beq\label{curvature} \beg{split}\Ric'&= \Ric -(d-2)(\Hess_{\log
f^{-1}}-(\d\log f^{-1})\otimes (\d\log f^{-1}))\\
&\qquad\qquad -(\DD \log f^{-1}+
(d-2)|\log f|_g^2)g\\
&= \Ric +(d-2)f^{-1}\Hess_f +(f^{-1}\DD f -(d-3)|\nn \log
f|_g)g.\end{split}\end{equation} Due to (\ref{N*1}) and
(\ref{curvature}), we are able to prove the following main result
in this subsection.

\beg{prp}\label{PP} For $Z\in C^2$, the Riemannian manifold
$(M',g')$ and the drift $Z':= f^2 Z+ \ff {d-2} 2 \nn f^2$ satisfies
$(H).$
\end{prp}

This Proposition will be implied by Lemma \ref{L*2} and Lemma
\ref{L*3} below. To prove these lemmas, we first clarify the
relationship between $\scr X_b^1(M',g)$ and $\scr X_b^1(M',g').$

\beg{lem}\label{L*1} For any $X\in \scr X(M'),$

\beq\label{N*2} \big| |\nn X|_g - |\nn' X|_{g'}\big| \le 3 |\nn f|_g
|X|_{g'}.\end{equation} Consequently,

\beq\label{N*3} f \scr X_b^1(M',g):=\{fX: \ X\in \scr
X_b^1(M',g)\}\subset \scr X_b^1(M',g')\subset \scr
X_b^1(M',g).\end{equation}
\end{lem}

\beg{proof}  For any $Y\in TM'$ with $|Y|_g=1$, one has $|f
Y|_{g'}=1$ and by (\ref{N*1}),

\beg{equation*}\beg{split} \big| |\nn_Y X|_g - |\nn'_{fY}
X|_{g'}\big|&= \big| |\nn_Y X|_g - |\nn'_{Y} X|_{g}\big|\le |\nn_Y
X- \nn'_YX|_g\\
&\le 3 |f^{-1} X|_g |\nn f|_g |Y|= 3 |X|_{g'}|\nn
f|_g.\end{split}\end{equation*} Thus, (\ref{N*2}) holds. Since $f$
is smooth with $0\le f\le 1,$ it is obvious that

$$\scr X_b^0(M',g')= f\scr X_b^0(M,g) \subset \scr X_b^0(M,g),$$
where $\scr X_b^0$ denotes the set of all bounded continuous vector
fields.

If $X\in \scr X_b^1(M',g),$ then (\ref{N*2}) implies

$$|\nn' (fX)|_{g'} \le |\nn (fX)|_g + 3|\nn f|_g|X|_g\le 4|\nn f|_g
|X|_g + f|\nn X|_g,$$ which is bounded. So, $f\scr X_b^1(M,g)\subset
\scr X_b^1(M',g').$

On the other hand, if $X\in \scr X_b^1(M',g'),$ then by (\ref{N*2}),

$$|\nn X|_g \le |\nn' X|_{g'} + 3 |\nn f|_g |X|_{g'}$$ is bounded.
Therefore, the proof is finished.\end{proof}

\beg{lem}\label{L*2} For any $C^2$ vector field $Z$ on $M$, one has
 $Z'\in \scr
X_b^2(M',g').$\end{lem}

\beg{proof} We shall prove $f^2 Z\in \scr X_b^2(M',g')$ and $\nn f^2
\in \scr X_b^2(M',g')$ respectively.

(a)  For any $X\in \scr X_b^1(M',g'),$ by (\ref{N*1}) we have

\beg{equation*}\beg{split} \nn'_X(f^2 Z) &=\nn_X (f^2 Z) -\<X,\nn
f\>_g f Z - \<Z, \nn f\>_g f X +\<fZ, X\>_g\nn f\\
&= f\{f \nn_XZ + \<Z,X\>_g\nn f +\<X,\nn f\>_g Z -\<Z, \nn f\>_g X\}
=: f U.\end{split}\end{equation*} By (\ref{N*3}) we have $X\in \scr
X_b^1(M',g)$. Moreover, $Z\in \scr X_b^2(M',g).$ Thus, $U\in \scr
X_b^1(M',g).$ So, by (\ref{N*3}), $\nn'_X (f^2 Z)\in \scr
X_b^1(M',g')$ for all $X\in \scr X_b^1(M',g').$ This means that $f^2
Z\in \scr X_b^2(M',g').$

(b) Let $X\in \scr X_b^1(M',g')$, it remains to prove that
$|\nn'\nn'_X\nn f^2|_{g'}$ is bounded. By (\ref{N*1})

\beg{equation*}\beg{split} \nn'_X\nn f^2 &= \nn_X\nn f^2 -2\<X, \nn
f\>_g \nn f- 2 |\nn f|_g^2 X + 2 \<\nn f, X\>_g \nn f \\
&= 2\<X,\nn f\>_g\nn f + 2 f \nn_X \nn f -2|\nn f|^2
X.\end{split}\end{equation*} By (\ref{N*3}), $f\nn_X\nn f\in \scr
X_b^1(M',g')$. So, it suffices to prove that

$$I:= |\nn' (\<X,\nn f\>_g \nn f - |\nn f|_g^2 X)|_{g'}$$
is bounded. By (\ref{N*2}),

\beg{equation*}\beg{split} I &\le |\nn (\<X,\nn f\>_g \nn f)|_g +
|\nn (|\nn f|_g^2X)|_g + 3 |\nn f|_g^2 |\<f^{-1}X,\nn f\>_g|+ 3|\nn
f|_g^3|X|_{g'}\\
&\le 5 |\nn^2 f|_g |\nn f|_g |X|_g + 5 |\nn f|_g^3 |X|_{g'} + 2|\nn
f|_g^2 |\nn X|_g\end{split}\end{equation*} which is bounded since
$X\in \scr X_b^1(M',g')\subset \scr X_b^1(M',g)$. \end{proof}

\beg{lem}\label{L*3} The Ricci curvature  $\Ric' \in C_b^1(M',g')$
and the curvature tensor $\Omega'\in C_b^0(M',g')$.\end{lem}

\beg{proof} By (\ref{curvature}), there exists a smooth $2$-tensor
$\scr T$ on $M$ such that

$$\Ric'(X,Y)= f\scr T(f^{-1}X, f^{-1}Y),\ \ \ \ X, Y\in \scr X_b^1(M',g').$$ Then $R'$
is bounded since $|\cdot|_g= f|\cdot|_{g'}$. Assuming
$|X|_{g'}=|Y|_{g'}=1,$ we obtain from the above formula and
(\ref{N*2}) that

\beg{equation*}\beg{split} &|\nn'\, \Ric' (X,Y)|_{g'} = f
|\nn\, \Ric' (X,Y)|_g \\
&\le f |\nn f|_g  |\scr T|_g
  + f |\nn \scr T|_g( |\nn (f^{-1}X)|_g+|\nn (f^{-1}Y)|_g)\\
&\le f |\nn f|_g   |\scr T|_g    +   |\nn \scr T|_g  \big(|\nn
X|_{g}+ |\nn Y|_g + 2 |\nn f|_{g}\big),\end{split}\end{equation*}
which is bounded on $M'$ since $X, Y\in \scr X_b^1(M',g')\subset
\scr X_b^1(M',g)$, $\scr T$ is smooth on $M$, and $M'\subset M$ is
relatively compact. Therefore $\Ric' \in C_b^1(M',g')$. The same
argument does work for $\Omega'$. The proof is finished.
\end{proof}

\subsection{Proof of Theorem \ref{T1.1}}

By Greene-Wu's approximation theorem \cite{GW}, we take a smooth
positive function $\tt\rr$ on $M$ such that

\beq\label{L1} |\tt\rr- \rr_o|\le 1, \ \ \ff 1 2\le |\nn \tt\rr| \le
2\ \ \text{and}\ (\DD+Z)\tt\rr\le (\DD+Z)\rr_o+1,\end{equation}
where the last inequality is restricted outside
$\{o\}\cup\text{cut}(o).$ Moreover, by the approximation theorem, we
may and do assume that $Z\in C^2$.

\beg{lem} \label{LL}   $(\ref{4.1})$ implies

$$(\DD+Z)\tt \rr  \le K+2+\psi(\tt\rr +1),\ \ \ \rr\ge 1.$$ \end{lem}

\beg{proof} For $x\notin \text{cut}(o)$ with $\rr_o=\rr_o(x)\ge 1$,
let $\ell : [0,\rr_o]\to M$ be the minimal geodesic from $o$ to $x$.
Let $U=\dot \ell$ and $\{U_i: 1\le i\le d-1\}$ be constant vector
fields along $\ell$ such that $\{U, U_i: 1\le i\le d-1\}$ is an
orthonormal basis. Let $J_i$ be the Jacobi field along $\ell$ with
$J_i(0)=0$ and $J_i(\rr_o)=U_i, 1\le i\le d.$ Let $h(s)=
1-(\rr_o-s)^+.$ By the second variational formula and the index
lemma, we have

\beq\label{L2}\beg{split} \DD\rr_o(x)&=
\sum_{i=1}^{d-1}\int_0^{\rr_o} \big\{|\nn_U J_i|_g^2 -R(J_i, U,
J_i,U)\big\}\\
&\le \sum_{i=1}^{d-1}\int_0^{\rr_o} \big\{|\nn_U (hU_i)|_g^2
-h^2R(U_i, U,  U_i, U)\big\} = 1 -\int_0^{\rr_o} h^2
\Ric(U,U).\end{split}\end{equation} Next,

$$Z \rr_o= \<Z,U\>_g(x)= \int_0^{\rr_o}\ff{\d}{\d s} \{h^2
\<Z,U\>_g\}\d s \le \int_0^{\rr_o} h^2\<\nn_U Z,U\>_g
+\psi\circ\rr_o(x).$$ Combining this with (\ref{L2}) we obtain

$$(\DD+Z)\rr_o\le K +1 +\psi\circ\rr_o.$$ Therefore, the proof is finished by
(\ref{L1}).\end{proof}

\ \newline \emph{Proof of Theorem \ref{T1.1}.} Let $h_0\in
C_b^\infty$ be decreasing such that $0\le h_0\le 1, h_0(s)=1$ for
$s\le 1$, and $h_0(s)=0$ for $s\ge 2$. Let

$$h_n(s):= h_0\bigg(\ff 1 n \int_0^s\ff {\d t} {\psi(t+1)}\bigg),\ \ \ s\ge 0, n\ge 2.$$
 For any $n\ge 2,$ let $f_n= h_n (\tt\rr).$ Since $\psi>0$ is smooth with
$\int_0^\infty \psi(s)^{-1}\d s=\infty,$ we have $f_n\in
C_0^\infty(M).$ Let $\mu_n$ be the distribution on $W_o(M)$ for the
diffusion process generated by $f_n^2L.$  Then $\mu_n\to \mu$
strongly; that is, for any bounded measurable function $F$ on
$W_o(M),$

\beq\label{L3} \lim_{n\to\infty}\mu_n(F)=\mu(F).\end{equation}
Indeed, letting $\tau_n$ be the hitting time of the $L$-diffusion
process to the set $\big\{\int_0^{\tt\rr} \psi(s)^{-1}\ge n \big\}$,
these two diffusion processes have the same distribution up to
$\tau_n.$ So,

$$|\mu(F)-\mu_n(F)|\le 2\|F\|_\infty \P(\tau_n\le T).$$ Since
 $\tau_n\to\infty$ as $n\to\infty,$ we obtain
(\ref{L3}). Then, it is standard that

\beq\label{L10} W_{2,d_\infty}^2(F\mu, \mu)\le
\liminf_{n\to\infty}W_{2,d_\infty}^2(F_n\mu_n, \mu_n)\end{equation}
for $F_n:= F/\mu_n(F).$

Now, let $\Ric^n, \nn^n$ be the Ricci curvature and Levi-Civita
connection induced by $g_n:= f_n^{-2}g$ on $M_n:= \{f_n>0\}.$ Let
$Z_n= f_n^2 Z+ (d-2)f_n\nn f_n.$ By (\ref{L10}), Propositions
\ref{P2} and \ref{PP}, it remains to prove

\beq\label{FF} \Ric^n -\nn^n Z_n\ge -K_n  \end{equation} for some
positive constants $K_n\to K$ as $n\to\infty.$ Let $X\in TM_n$ with
$|X|_{g_n}=1.$  By (\ref{curvature}), we have

$$\Ric^n(f_nU,f_nU)= f_n^2\Ric(U,U) +(d-2)f_n\Hess_f(U,U) +f_n\DD
f_n -(d-3)|\nn f_n|^2,\ \ \ |U|_g=1.$$ Combining this with  the
first display on \cite[page 114]{TW}, we obtain

\beg{equation*}\beg{split} &(\Ric^n -\nn^n Z_n)(f_n U, f_n U)\\
& \ge f_n^2 (\Ric -\nn Z)(U,U) + f_n (\DD+Z)f_n -c_1(|\nn f_n|_g^2+
|Z|_g |\nn f_n|_g),\ \ \ |U|_g=1\end{split}\end{equation*} for some
constant $c_1>0.$  Combining this with (\ref{4.1}) we obtain

$$\Ric^n -\nn^n Z_n  \ge -K + f_n
(\DD+Z)f_n -c_1(|\nn f_n|_g^2+ |Z|_g |\nn f_n|_g).$$ Therefore, to
ensure (\ref{FF}) it suffices to show that

\beq\label{Last} \lim_{n\to\infty}\inf\{f_n (\DD+Z)f_n -c_1(|\nn
f_n|_g^2+ |Z|_g |\nn f_n|_g)\}=0.\end{equation}  By Lemma \ref{LL},
$h_0'\le 0$ and  $|\nn \tt\rr|\le 2$,

$$(\DD+Z)f_n \ge  - \ff{\|h'_0\|_\infty (K+2+\psi(\tt\rr
+1))}{n\psi(\tt\rr+1)} -\ff{2\|h_0''\|_\infty}{n^2\psi(\tt\rr +1)^2}
$$  which goes to zero uniformly as $n\to\infty$.  Similarly,
$|\nn f_n|_g^2 +|Z|_g|\nn f_n|_g\to 0$ uniformly too. \qed

\section{An extension to free path spaces}

Let $\nu$ be a probability measure on $M$ such that

\beq\label{5.1} W_{2,\rr}(f\nu,\nu)^2\le C_0\nu(f\log f),\ \ \ f\ge
0, \nu(f)=1\end{equation} holds for some constant $C_0>0.$ Let
$P_\nu$ be the distribution of the $L$-diffusion process starting
from $\nu$ up to time $T>0$, which is then a probability measure on
the free path space $W(M)=C([0,T];M).$

\beg{thm}\label{T5.1} Under $(\ref{4.1})$ and the growth condition
for $|Z|$ stated in Theorem $\ref{T1.1}$ for some $($and hence
any$)$ fixed point $o\in M.$ Then

$$W_{2,d_\infty}(FP_\nu, P_\nu)^2\le \Big(C_0 \e^{KT}+
2\ff{\e^{KT}-1}K\Big)P_\nu(F\log F),\ \ \ \ F\ge 0,
P_\nu(F)=1.$$\end{thm}

\beg{proof} (a) Without loss of generality, we assume that $F\in
\FC$ is strictly positive. Let $P_x$ be the distribution of the
$L$-diffusion process starting from $x$, and let $f(x)=P_x(F), F_x=
\ff F{f(x)}.$ Then $\nu(f)=P_x(F_x)=1$ and

\beq\label{5.2} P_{f\nu}= \int_M (F_x P_x)f(x)\nu(\d x),\ \ \
P_{f\nu}=\int_MP_x f(x)\nu(\d x).\end{equation}. By the triangle
inequality,

\beq\label{5.3} W_{2,d_\infty}(FP_\nu, P_\nu)\le
W_{2,d_\infty}(FP_\nu, P_{f\nu})+ W_{2,d_\infty}(P_{f\nu},
P_\nu).\end{equation}

(b) It is well-known that in a class of probability measures on a
Polish space with bounded second moment, the weak convergence is
equivalent to the convergence in the $L^2$ Wasserstein distance
(see e.g. \cite{Ra}). Noting that $x\mapsto P_x$ and $x\mapsto
F_xP_x$ are continuous in the weak topology for probability
measures on $W(M)$, and due to (\ref{4.1}), $\sup_x
P_x(\e^{\d_\infty(x,\cdot)})<\infty,$ we conclude that

$$x\mapsto W_{2,d_\infty}(P_x, F_xP_x)$$
is continuous. Furthermore, Theorem \ref{T1.1} and the uniform
boundedness of $F_x$ imply that this function is bounded. Therefore,
it is is to see from (\ref{5.2}) that

\beq\label{5.4} W_{2,d_\infty}(FP_\nu, P_{f\nu})^2\le
\int_MW_{2,d_\infty}(F_xP_x, P_x)^2f(x)\nu(\d x).\end{equation}
Indeed, letting $\{A_{i,n}:\ i\ge 1\}_{n\ge 1}$ be a sequence of
measurable partitions of $M$ such that

$$\nu(A_{i,n})+ \text{dia}(A_{i,n})\le \ff 1 n, \ \ \ i,n\ge 1,$$
where $\text{dia}(A_{i,n})$ is the diameter of $A_{i,n}.$ By the
continuity of $f$, let $x_{i,n}\in \bar A_{i,n}$ such that

$$f(x_{i,n})\nu(A_{i,n})=\int_{A_{i,n}}f(x)\nu(\d x),\ \ \ \ i,n\ge
1.$$ Let $\pi_{i,n}\in \scr C(F_{x_{i,n}}P_{x_{i,n}}, P_{x_{i,n}})$
such that

$$\int_{W(M)\times W(M)} d_\infty^2\d\pi_{i,n}= W_{2, d_\infty}(F_{x_{i,n}}P_{x_{i,n}},
P_{x_{i,n}})^2,\ \ \ i,n\ge 1.$$ Then

$$\pi_n:= \sum_{i=1}^\infty f(x_{i,n})\nu(A_{i,n})\pi_{i,n}\in \scr C((FP_\nu)_n,
(P_{f\nu})_n),$$ where

$$(FP_{\nu})_n:= \sum_{i=1}^\infty f(x_{i,n}) \nu(A_{i,n})
F_{x_{i,n}} P_{x_{i,n}}\to FP_\nu$$ and

$$(P_{f\nu})_n:=\sum_{i=1}^\infty f(x_{i,n}) \nu(A_{i,n})
 P_{x_{i,n}}\to P_{f\nu}$$  weakly as $n\to 0$, then

 \beg{equation*}\beg{split}
 W_{2,d_\infty}(FP_\nu, P_{f\nu})^2&= \lim_{n\to\infty}W_{2,d_\infty}((FP_\nu)_n,
 (P_{f\nu})_n)^2 \\
 &\le\lim_{n\to\infty} \int_{W(M)\times W(M)} d_\infty^2 \d\pi_n\\
 & =
  \lim_{n\to\infty}\sum_{i=1}^\infty
 f(x_{i,n})\nu(A_{i,n})W_{2,d_\infty}(F_{x_{i,n}}P_{x_{i,n}},
 P_{x_{i,n}})^2\\
 & =\int_MW_{2,d_\infty}(F_xP_x, P_x)^2f(x)\nu(\d x),\end{split}\end{equation*}
 Therefore, (\ref{5.4}) holds.  Combining this with Theorem \ref{T1.1}, we obtain

 \beq\label{5.5} \beg{split}W_{2,d_\infty}(FP_\nu, P_{f\nu})^2 &\le
 \ff{2(\e^{KT}-1}K \int_{W(M)\times M} \{F_x(\gg)\log
 F_x(\gg)\}P_x(\d\gg) f(x)\nu(\d x)\\
 &=\ff{2(\e^{KT}-1}K
 \big(P_\nu(F\log F)- \nu(f\log f)\big).\end{split}\end{equation}

 (c) To estimate $W_{2,d_\infty}(P_{f\nu}, P_\nu),$ let $\hat\pi\in \scr (f\nu,\nu)$ such that

$$W_{2,\rr}(f\nu,\nu)^2=\int_{M\times M}\rr^2\d\hat\pi,$$ and let
$(X_t,Y_t)$
 be the coupling by parallel displacement for the $L$-diffusion
 process with initial distribution $\hat \pi.$ By
 \cite[(3.2)]{PTRF} (note that the present $L$ is half of the one
 therein)

 $$\rr(X_t, Y_t)\le \rr(X_0,Y_0) \e^{Kt/2},\ \ \ t>0.$$ Thus,

 $$W_{2,d_\infty}(P_{f\nu}, P_\nu)^2\le \mathbb E\max_{t\in [0,T]}
 \rr(X_t,Y_t)^2 \le \e^{KT} \mathbb E \rr(X_0,Y_0)^2 = \e^{KT}
 W_{2,\rr}(f\nu,\nu)^2.$$
 Then it follows from (\ref{5.1}) that

 $$W_{2,d_\infty}(P_{f\nu}, P_\nu)^2\le C_0\e^{KT}\nu(f\log f).$$
 Combining this with (\ref{5.3}) and (\ref{5.5}) we arrive at

 \beg{equation*}\beg{split} &W_{2,d_\infty}(FP_\nu, P_\nu)^2 \le
 (1+\dd) W_{2,d_\infty}(FP_\nu, P_{f\nu})^2 + (1+\dd^{-1}) W_{2,
 d_\infty}(P_{f\nu}, P_\nu)^2\\
 & \le \ff {2(1+\dd)(\e^{KT}-1)}{K} P_\nu(F\log F) +\Big( C_0
 (1+\dd^{-1}) \e^{KT} - \ff{2(1+\dd) (\e^{KT}-1)}K\Big)\nu(f\log
 f).\end{split}\end{equation*}
 Then the proof if finished by taking  $\dd= C_0
 K\e^{KT}/2(\e^{KT}-1).$ \end{proof}

\section{Appendix: regularity of $A^I$}

\medskip
Let $V$ be a smooth manifold. For the convenience of our
exposition, we shall introduce $V$-valued smooth Wiener functional
in the following way (for a general definition, we refer to
\cite{Ma}, p.78).
\medskip
Let $\Phi: W_o(\R^d)\rightarrow V$ be a measurable map. Let $p>1$.
We say that $\Phi$ is derivable if there exists $\nabla\Phi(w)\in
\H\otimes T_{\Phi(w)}V$ satisfying $\E(||\nabla\Phi||_{\H\otimes
TV}^p)<+\infty$ such that for each $h\in\H$, $\Phi$ admits a
version $\Phi_h$ such that $\varepsilon\ra \Phi_h(w+\varepsilon
h)$ is $C^1$ and
$$ \ff{\d}{\d \vv} \Phi_h(w +\vv
h)|_{\vv=0}=\nabla\Phi(w)\cdot h\in T_{\phi(w)}V.$$ Then
$\nabla\Phi$ is a map from $W_o(\R^d)$ into $\H\otimes TV$.
Inductively, we define high order derivatives $\nabla^k\Phi:
W_o(\R^d)\rightarrow \H^{\otimes k}\otimes TV$. We say that
$\Phi\in \mathbb{D}_k^\infty$ if $\E(||\nabla^r\Phi||^p)<+\infty$
for all $r\leq k$ and $p>1$. We say that $\Phi$ is non-degenerated
in Malliavin sense if $\hbox{\rm
det}^{-1}[\nabla\Phi(\nabla\Phi)^*]\in \cap_{p>1}L^p$, where
$(\nabla\Phi(w))^*: T_{\Phi(w)}V\rightarrow \H$ is defined by
$$\bigl<(\nabla\Phi(w))^*v,h\bigr>_H=\bigl<\nabla\Phi(w)h,v\bigr>_{T_{\Phi(w)}V}.$$
The following result holds (see \cite{Ma}, chapter III).

\begin{thm}\label{thm5.1}
Let $\Phi\in\mathbb{D}_2^\infty$ be a $V$-valued non-degenerated
Wiener functional and $G\in
\mathbb{D}_1^\infty(W_o(\R^d),\mathbb{R})$, then the conditional
expectation $z\ra \E(G|\Phi=z)$ admits a continuous version.
\end{thm}

\medskip
Now we are going to prove the regularity of $A^I$.
\medskip

\begin{lem}\label{lem4.1}
Assume $(H)$. Then  the It\^o functional $u_t: W_0(\R^d)\to O(M)$
defined by $(1.2)$ belongs to $\mathbb{D}_1^\infty$.
\end{lem}

\begin{proof} We first note that for any $h\in \mathbb H$, the law
of $w\mapsto u_t(w+\vv h)$ is equivalent to that of $u_t$ and
furthermore,

\begin{equation*}
\beta(t):=\bigl<\theta, D_h u_t\bigr>,\ \ \ \rho(t):=\bigl<\Theta,
D_h u_t\bigr>
\end{equation*} satisfy (see \cite[(2.21)]{Bi})

\begin{equation}\label{4.10}
\beg{cases} \d\beta(t) =(\dot h(t)+ \{\nabla Z\}^\#(u_t)\,\beta(t)
-\rho(t) Z^\#(u_t))\d t + \rho(t)
(Z^\#(u_t) \d t + \circ \d w_t),\\
\d\rho(t) =\OO_{u_t}(u_t^{-1}Z_{\pi u_t} \d t + \circ \d w_t,
\beta(t)).
\end{cases}
\end{equation}
Here, $Z^\#(u):=\<Z, u\cdot\>\in \R^d$ for $u\in O(M)$,
$(\theta,\Theta)$ is the parallelism of $O(M)$, an
$\mathbb{R}^d\times \hbox{\rm so}(d)$-valued one-form on $O(M)$
defined by

$$\theta_u(\tt X)= u^{-1} \pi^* \tt X,\ \ \ \Theta_u(\tt X)= q_u^{-1}(P_V \tt X),\ \ \ u\in O(M), \tt X\in T_u O(M),$$
where $P_V$ is the orthogonal projection from $TO(M)$ onto the space
of vertical tangent vectors on $O(M)$, and

$$q_u:\ \text{so}(d)\ni \aa\mapsto \ff{\d}{\d s} \{u \e^{-s\aa}\}|_{s=0}\in P_V T_u O(M)$$ is an endomorphism.

Let $D(d)=\R^d\times \hbox{so}(d)$. For $r\in O(M)$, we denote by
${\mathcal M}_j(u)$ the endomorphism of $D(d)$ defined by
$$  (x,A)\ra \bigl(\{\nabla Z\}^\#(u)\cdot x+Ae_j,
\OO_r(e_j,x)\bigr).$$    Let $J_{t,s}$ solve the equation on $\scr
L(\R^d; \R^d\times \text{so}(d)):$

\begin{equation}\label{4.11}
\frac{d}{\d t}J_{t,s}=\Bigl(\sum_{j=1}^d {\mathcal M}_j(u_t)\circ
\bigl[(u_t^{-1} Z_{\pi u_t})^j \d t+\d
w_t^j\bigr]\Bigr)J_{t,s},\quad t>s, \ J_{s,s}=(\Id_{D(d)}, 0).
\end{equation}
Then  (see \cite[page 292]{Ma})
$$
(\beta(t),\rho(t))=\int_0^t J_{t,s}\dot h(s)\,\d s.
$$
This completes the proof due to the fact that

$$\Big|\ff{\d}{\d \vv} u_t(w +\vv h)\big|_{\vv=0}\Big|^2_{T_{u_t}O(M)}
=|\bb(t)|^2 + |\rr(t)|^2$$ and the boundedness of $\mathcal M, \OO$
and $Z$.
\end{proof}

\begin{lem}\label{lem4.2} Assume $(H)$.
Let $Q_t= Q_{t,0}.$ Then

\begin{equation}\label{(4.12)}
D_hQ_t=-\frac{1}{2}\int_0^t Q_{t,s}\{\nn_{\pi^*
D_hu_s}\Ric_Z\}^\#(u_s)\, Q_s\, \d s,\ \ \ h\in\H.
\end{equation} Consequently, $Q_t\in \mathbb{D}_1^p(W_0(\R^d))$ for $p\ge 1.$
\end{lem}

\begin{proof}
Differentiating (\ref{4.4}),  we obtain
\begin{equation*}
\frac{\d D_h Q_t}{\d t}=-\frac{1}{2}\{\nn_{\pi^*
D_hu_t}\Ric_Z\}^\#(u_t) Q_t-\frac{1}{2}(\Ric_Z^\#(u_t)) D_h Q_t,\ \
D_h Q_0=0.
\end{equation*}
So, we get the expression (\ref{(4.12)}) and thus, $Q_t\in
\mathbb{D}_1^p(W_0(\R^d))$ for $p\ge 1$ due to $(H)$ and Lemma
\ref{lem4.1}.
\end{proof}

\begin{prp}\label{lem4.3} Assume $(H).$ Then $A^I: T M^I\to TM^I$
has a
 continuous $\mu_I$-version.
\end{prp}

\begin{proof} Let $K_{ij}(s)=Q_{s_j,s}Q_{s_i,s}^*$.
Note that $u_{s_i}^{-1}X_i(\gamma(s_i))=\theta(X_i^\#)_{u_{s_i}}$.
Then, for any compactly supported smooth vector fields $X, Y$ on
$M^I$,

$$ \bigl<u_{s_j}Q_{s_j,s}Q_{s_i,s}^*u_{s_i}^{-1}X_i ,\ Y_j \bigr>_g
=\bigl<K_{ij}(s)\theta(X_i^\#)_{u_{s_i}},
\theta(X_j^\#)_{u_{s_j}}\bigr>:=G_{ij}(t).
$$

By lemmas \ref{lem4.1} and \ref{lem4.2}, $G_{ij}(t)$ are in
$\mathbb{D}_1^\infty(W_0(\R^d),\mathbb{R})$, so

$$ G:=\sum_{i,j=1}^N \int_0^{s_i\wedge s_j} G_{ij}(s)\,ds\in \mathbb{D}_1^\infty(W_0(\R^d),\mathbb{R}).$$
  By Theorem
\ref{thm5.1}, $z\mapsto \< A^I(z) X(z), Y(z)\>_{g^I}=
\int_{W_0(\R^d)} G\,\P(z,\d w)$ has a continuous version.
\end{proof}

\medskip
\begin{prp}\label{lem4.4}
Assume that $\Ric-\nn Z\leq K_1$, then $A^I$ is uniformly elliptic
with respect to $g^I$.
\end{prp}

\begin{proof}
Let $a=(a_1, \cdots, a_N)\in T_zM^I$. Suppose without losing the
generality, that $|a_N|=\max_{1\le i\le N} |a_i|$. Take $(X_1,
\cdots, X_N)$ be vector fields around $(z_1, \cdots, z_N)$ such that
$$(X_1(z_1), \cdots, X_N(z_N))=(a_1, \cdots, z_N).$$
We have

\begin{equation*}
\bigl<A^I(z)a,a\bigr>=\E_{\mu}\Bigl(\int_0^T|\sum_{j=1}^N
Q_{s_j,s}^*(u_{s_j}^{-1}X_j(\gamma(s_j)){\bf 1}_{(s<s_j)}|^2\,\d
s\Big|\Lambda_I=z\Bigr).
\end{equation*}

Let $s_{N-1}\leq s<s_N$ and $v\in\R^d$. Then by the assumption on
the upper bound of $\Ric$,
\begin{equation*}
\frac{d}{\d t}|Q_{t,s}v|^2\geq -K_1|Q_{t,s}v|^2.
\end{equation*}

It follows that $  |Q_{t,s}v|^2\geq \e^{-K_1(t-s)}|v|^2\geq
\e^{-K_1(s_N-s_{N-1})}|v|^2$. Therefore,

\begin{equation*}
\begin{split}
&\int_0^T|\sum_{j=1}^N
Q_{s_j,s}^*(u_{s_j}^{-1}X_j(\gamma(s_j)){\bf 1}_{(s<s_j)}|^2\,\d s\\
&\geq
\int_{s_{N-1}}^{s_N}|Q_{s_N,s}^*(u_{s_N}^{-1}X_N(\gamma(s_N))|^2\,
\d s\\
&\geq |X_N(\gamma(s_N))|^2 \e^{-K_1(s_N-s_{N-1})}(s_N-s_{N-1}).
\end{split}
\end{equation*}
Hence $$\bigl<A^I(z)a,a\bigr>\geq
|a_N|^2\e^{-K_1(s_N-s_{N-1})}(s_N-s_{N-1})\ge |a|^2
N^{-1}\e^{-K_1(s_N-s_{N-1})}(s_N-s_{N-1}).$$
\end{proof}

\paragraph{Acknowledgements}
The authors would like to thank the referee for useful comments on
an earlier version of the paper.

\end{document}